\documentclass[12pt]{article}

\usepackage{amsmath}
\usepackage{multirow} 
\usepackage{latexsym,amsfonts}
\usepackage{graphicx}
\usepackage{subcaption}
\usepackage{pdflscape}
\usepackage{caption}
\usepackage{tikz-cd} 
\usepackage{tikz}
\usepackage{enumitem} 
\usepackage{authblk}
\newtheorem{thm}{Theorem}[section]
\newtheorem{note}{Note}
\newtheorem{lem}[thm]{Lemma}

\newcommand{\Proof}{\medskip\noindent{\bf Proof.}\quad}
\newcommand{\qed  }{\hfill$\Box$}

\author[1*]{Sanja Rukavina}  
\author[2]{Vladimir D. Tonchev}
\affil[1]{Faculty of Mathematics, University of Rijeka, 51000 Rijeka, Croatia\\

E-mail: sanjar@math.uniri.hr, ORCID: 0000-0003-3365-7925}
\affil[2]{ Department of Mathematical Sciences, Michigan Technological University, Houghton, MI 49931, USA\\

E-mail: tonchev@mtu.edu, ORCID: 0000-0003-1806-3571}
\title{Extremal ternary self-dual codes of length 36 and symmetric 2-$(36,15,6)$ designs with an
automorphism of order 2}

\date{}
\begin{document} 
\maketitle

* Corresponding author.\\

\begin{abstract} 
In this note we report the classification of all symmetric  2-$(36,15,6)$ designs that admit an
automorphism of order 2 and their incidence matrices generate an extremal ternary self-dual
code. It is shown that up to isomorphism, there exists only one  such design, 
having a full automorphism  group of order 24, and the ternary code spanned by
its incidence matrix is equivalent to the Pless symmetry code.

\end{abstract}

\vspace*{0.5cm}



\section{Introduction}

We assume familiarity with the basic facts and notions 
from error-correcting codes, combinatorial designs and Hadamard matrices
\cite{BJL}, \cite{BBFKKW}, \cite{HP}, \cite{IS}, \cite{ton88}.
All codes considered in this paper are ternary.

The minimum distance, or equivalently, the minimum weight $d$
of a ternary self-dual code of length $n$ divisible by 12
satisfies the upper bound $d\le n/4 + 3$ (cf., e.g. \cite[9.3]{HP}),
and a self-dual code with minimum distance $d=n/4+3$ is called 
{\it extremal}.  The extremal self-dual codes support combinatorial 5-designs
by the Assmus-Mattson theorem \cite{AM},  \cite[8.4]{HP}.

Ternary extremal self-dual codes are known for the following lengths $n\equiv 0 \pmod {12}$:
$n=12$: the extended Golay code, unique up to equivalence; $n=24$: there are exactly two inequivalent codes,
the extended quadratic-residue code \cite{AM} and the Pless symmetry code $C(11)$ \cite{Pless69}, \cite{Pless72};
$n=36$: only one code is known, namely the Pless symmetry code $C(17)$ \cite{Pless69}, \cite{Pless72};
$n=48$: two codes are known, the extended quadratic-residue code and the Pless symmetry code $C(23)$;
$n=60$: three codes are known: the extended quadratic-residue code, the Pless symmetry code $C(29)$,
and a code found by Nebe and Villar \cite{NV}.

Huffman \cite{Huf} proved that any extremal ternary self-dual code of length 36
that admits an automorphism of prime order $p>3$
is monomially equivalent to the Pless symmetry code. More recently, Eisenbarth and Nebe \cite{EN}
extended Huffman's result by proving that the Pless symmetry code is the unique (up to monomial equivalence)
ternary extremal self-dual  code  of length 36 that admits an automorphism of order 3.
In addition, it was proved in \cite[Theorem 5.1]{EN} that if $C$ is an extremal ternary self-dual code
of length 36 then either $C$ is equivalent to the Pless symmetry code or the full automorphism group
of $C$ is a subgroup of the cyclic group of order 8. 

It is known \cite{Pless72} that the Pless symmetry code $C(q)$ of length $n=2q+2$, where $q\equiv -1 \pmod 3$ is an odd prime power, contains a set of $n$ codewords of weight $n$, which after replacing every entry equal to 2 by $-1$ form the rows
of a Hadamard matrix equivalent to the Paley-Hadamard matrix of type II.
 In particular, the Pless symmetry code $C(17)$ contains the rows of a Hadamard matrix $P$ of  Paley type II, 
 having a full automorphism group of order  $4\cdot 17(17^2 -1)=19584$, and the rows of $P$  span the code $C(17)$. 

It was shown in \cite{Ton21} that  the code $C(17)$ contains a second equivalence class of Hadamard 
matrices of order 36 having as rows codewords of $C(17)$. Any matrix $H$ from the second equivalence
class has a full automorphism group of order 72 and the rows of $H$ span the code $C(17)$.
In addition,  $H$ is monomially equivalent to a regular  Hadamard matrix $H'$,
such that every row of $H'$ has 15 entries equal to 1 and 21 entries equal to $-1$ .
The symmetric 2-$(36,15,6)$ design $D$ with $(0,1)$-incidence obtained by replacing every entry $-1$ of $H'$ 
with 0 has a trivial full automorphism group, and the row span of its incidence matrix over $GF(3)$ is
equivalent to the Pless symmetry code $C(17)$.

In Section 2, we consider codes that are obtained from $C(17)$  by negating 
code coordinates so that the resulting new code contains the all-one vector,
and examine Hadamard matrices formed by codewords of weight 36.
This examination reveals the surprising fact that the Pless symmetry code $C(17)$ is equivalent to
a code  spanned by the rows of a regular Hadamard matrix which
is monomially equivalent to the Paley-Hadamard matrix of type II,  the associate
symmetric 2-$(36,15,6)$ design has a full automorphism group of order 24,
and its incidence matrix spans a code equivalent to $C(17)$.

Motivated by this phenomenon and the results from  \cite{EN}, in Section 3 we classify
all symmetric 2-$(36,15,6)$ designs that admit an automorphism of order 2
and their incidence matrices span an extremal ternary self-dual code of length 36.
The results of this classification show that, up to isomorphism, there is only one such
design, being isomorphic to the design associated with the regular presentation of the
Paley-Hadamard matrix of type II found in Section 2.

\section{Hadamard matrices derived from the symmetry code of length 36}
\label{sec2}

In this section we consider Hadamard matrices of order 36 whose rows are obtained from
codewords of full weight 36  in the ternary Pless symmetry code  $C(17)$
after replacing all codeword entries that are equal to 2 with $-1$. 
The code $C(17)$ contains exactly 888 codewords of weight 36. 
In what follows, we denote by $U$ the set of all  888 codewords in $C(17)$ of full weight.
 The 3-rank (that is, the rank over the finite field of order 3, $GF(3)$) of the matrix having as rows the
 codewords from $U$ is 18, hence $U$  spans  $C(17)$.
In the context of Hadamard matrices, we view the codewords from $U$
as vectors with components $\pm 1$. 
 
 \begin{lem}
 \label{l1} 
 A set  $S$ of 36 codewords from $U$
 is the row set of a Hadamard matrix of order 36  if and only if
 the Hamming distance between every two codewords from $S$ is 18.
 \end{lem}
 \Proof The inner product of two vectors 
 $u=(u_1,\ldots, u_{36})$, $v=(v_1, \ldots , v_{36})$ of length 36 with entries $\pm 1$ is zero
 over the field of real numbers if and only if there are exactly 18 indices $i$, $1\le i\le 36$,  such that
 $u_{i}v_{i}=-1$, hence $u_i \neq v_i$, and  the remaining 18 indices $j$, $1\le j\le 36$, satisfy $u_{j}v_{j}=1$,
 that is, $u_j = v_j$.
 \qed
 
 By Lemma \ref{l1}, every Hadamard matrix of order 36 whose rows are codewords from $U$
 corresponds to a clique of size 36 in a graph $\Gamma$ having as vertices the codewords from $U$,
where two codewords are adjacent in $\Gamma$ if they are at Hamming distance 18 from each other.
Enumerating all Hadamard matrices up to equivalence by a direct search for cliques of size 36 is computationally
 infeasible because for every Hadamard matrix $H$, the graph $\Gamma$ contains $2^{36}$ distinct 36-cliques 
 that correspond to  $2^{36}$ Hadamard matrices that are
 equivalent to $H$ and are obtained by negating rows of $H$. Note that
 negating a row  is equivalent to a scalar multiplication of the corresponding codeword by 2 (mod 3), 
thus any negation of a row preserves $U$.

A simple way to reduce the search is by considering normalized Hadamard matrices.

\begin{lem}
\label{l2}
If $H$ is a Hadamard matrix whose rows are codewords from $U$
then the row set of any normalized Hadamard matrix obtained from $H$ 
belongs to a code which is monomially equivalent to $C(17)$.
\end{lem}
\Proof The matrix $H$ is normalized with respect to  column $j$
by negating (or multiplying by 2 (mod 3)) all rows of $H$ having entry $-1$ in the $j$th column.
Any such negation preserves $U$. Similarly,  $H$ is normalized with respect to row $i$
by negating all columns of $H$ having entry $-1$ in the $i$th row. These negations of
 the columns of $U$ either preserve $U$ or change $U$ to a new set that spans
a linear code which is monomially equivalent to $C(17)$.
\qed

The set $U$ of the Pless symmetry code $C(17)$, as defined in
\cite{Pless69}, \cite{Pless72}, does not contain the all-one codeword $\bar{1}=(1,\ldots,1)$,
 while it contains a codeword $v$ with one entry equal to  $-1$,
located in the code coordinate labeled by $\infty$, and 35 entries equal to 1.
Negating (that is, multiplying by 2 (mod 3)) the code coordinate of $C(17)$ labeled by $\infty$ 
transforms the symmetry code $C(17)$ into a monomially equivalent code
 $L(17)$ which does contain the all-one vector $\bar{1}$ \cite{Ton21}.
 If $x \in L(17)$ is a codeword of full weight 36, we denote by $w_{i}(x)$ the number 
 of entries of $x$ that are equal to $i$ ($i = 1, 2$).
 
 The complete weight distribution of the set $W$ consisting of all 888 codewords of $L(17)$ of weight 36  
is listed in Table \ref{tab1}. The set $W$ is available at 
\begin{verbatim}
 https://pages.mtu.edu/~tonchev/W.txt
 \end{verbatim}

As shown in \cite{Ton21}, the set $M$ of the 70 codewords $x\in W$ with weight structure $(w_{1}(x),w_{2}(x))=(18,18)$
 form the $(1,2)$-incidence
matrix of a Hadamard 3-$(36,18,8)$ design with full automorphism group of order $272=2^{4}17$, 
while adding to $M$ the all-one vector $\bar{1}$ and the all-two vector
$\bar{2}=2\bar{1}$ gives a set of 72 vectors that comprises of the rows of a normalized Hadamard matrix and its negative,
being monomially equivalent to the Paley-Hadamard matrix of type II, having full monomial automorphism group of order
$19584=2^{7}3^{2}17$.

In addition, it was  shown in  \cite{Ton21} that the set of 408 codewords $x\in W$  with weight structure
 $(w_{1}(x),w_{2}(x))=(15,21)$
contains exactly 272 subsets of size 36 that are regular Hadamard matrices of order 36.
All these 272 regular Hadamard matrices are pairwise equivalent, and have an automorphism group of
order 72, while the associated symmetric 2-$(36,15,6)$ design  has trivial full automorphism group,
and its $(0,1)$-incidence matrix has 3-rank 18, hence spans the code $L(17)$.

All inequivalent Hadamard matrices of order 36 whose rows are obtained from codewords in $L(17)$ of weight 36
can be enumerated by using Lemma \ref{l2} as follows. We consider $W$ as an $888 \times 36$ matrix, and
 for every integer $i$, $1\le i \le 888$, we define a  matrix $W_i$, being the matrix obtained by negating all columns
 of $W$ which have entry $-1$ in row $i$. Thus, the $i$th row of $W_i$ is the all-one vector $\bar{1}$,
 and this row must be a row of every normalized Hadamard matrix consisting of rows of $W_i$.
 We refer to $W_i$ as a matrix obtained by switching of $W$ with respect to row $i$.
  To reduce the search further, we consider only normalized Hadamard matrices with first column being the
  all-one column.

Next we define a graph $\Gamma_i$ having as vertices all rows of $W_i$ with first entry 1 and exactly 18 entries
equal to 1, where two vertices are adjacent in $\Gamma_i$ if and only if the Hamming distance between the
corresponding rows of $W_i$ is 18.

\begin{lem}
\label{l3}
(a) Any set of 35 rows of $W_i$ that corresponds to a clique of size 35 in $\Gamma_i$, 
together with the all-one row $\bar{1}$,
is the set of rows of a normalized Hadamard matrix.\\
(b) The maximum clique size in $\Gamma_i$ is 35.
\end{lem}

\Proof Part (a) follows from Lemma \ref{l1}.
If $H$ is a normalized Hadamard matrix of order $n=2d$, then deleting the all-one column from $H$
gives a $2d\times (2d-1)$ matrix whose rows form a code $C$ of length $2d-1$ and minimum distance $d$
over the alphabet $\{ 1, -1 \}$ that meets the Plotkin upper bound \cite[2.2]{HP}, \cite[2.1]{ton88}:
\[ |C| \le 2d. \]
This proves part (b).
\qed

\begin{table}
\begin{center}  
\begin{tabular}{|r|c|}
\hline
\#$x$ & $(w_{1}(x),w_{2}(x))$\\
\hline
\hline
1 & (0,36) \\
\hline
408 & (15,21) \\
\hline
70 & (18,18) \\
\hline
408 & (21,15) \\
\hline
1 & (36,0)\\
\hline
\end{tabular}
\caption{The complete weight distribution of $W$}\label{tab1}
\end{center}
\end{table}

One can compute representatives of the equivalence classes of Hadamard matrices
by  an examination of the 35-cliques in the graphs $\Gamma_i$, $1\le i \le 888$.
Since every codeword $x\in W$ with weight structure $(w_{1}(x),w_{2}(x))=(21,15)$
is the negative of the codeword $2x \in W$ with weight structure $(15,21)$,
it is sufficient to examine the graphs $\Gamma_i$ such that the weight structure of the $i$th row of $W$
is $(15,21)$ or $(18,18)$. The incidence structure with $(0,1)$-incidence matrix obtained from $W$
by replacing all entries that are equal to 2 with zero, has full automorphism group $G$ of order 272.
The group $G$ partitions the set of  408 rows of $W$ with weight structure $(15,21)$ into four orbits
of lengths 68, 68, 136, 136, and the set of 70 rows with weight structure $(18,18)$ into four orbits
of lengths 2, 17, 17, 34. Therefore, it is sufficient to examine eight switchings of $W,$
one for each orbit, all together.

An examination of the matrices  $W_i$ obtained
by switching of $W$ with respect to a row $i$ having weight structure $(18,18)$
shows that any such matrix has the same complete weight distribution as $W$ 
(see Table \ref{tab1}),
and the resulting unique (up to a permutations of rows) normalized Hadamard matrix
is equivalent to the Paley-Hadamard matrix of type II.
 
 An examination of the matrices  $W_i$ obtained
by switching of $W$ with respect to any of the 408 rows having weight structure $(15,21)$
shows that any such matrix has a complete weight distribution given in Table \ref{tab2},
where $W_{408}$ is obtained by the switching of W with respect to row no. 408.

\begin{table}
\begin{center}  
\begin{tabular}{|r|c|}
\hline
\#$x$ & $(w_{1}(x),w_{2}(x))$\\
\hline
\hline
1 & (0,36) \\
\hline
93 & (12,24) \\
\hline
36 & (15,21) \\
\hline
628 & (18,18) \\
\hline
36 & (21,15)\\
\hline
93 & (24,12)\\
\hline
1 & (36,0)\\
\hline
\end{tabular}
\caption{The complete weight distribution of $W_{408}$}\label{tab2}
\end{center}
\end{table}

Any graph $\Gamma_i$ associated 
with a matrix  $W_i$ obtained
by switching of $W$ with respect to a row $i$ with weight structure $(15,21)$
has 314 vertices and contains exactly 24 cliques of size 35. The resulting 24
normalized Hadamard matrices are all equivalent to the Hadamard matrix 
with full automorphism group of order 72 that was discovered in \cite{Ton21}
as a regular Hadamard matrix derived from the code $L(17)$,
 whose related symmetric 2-$(36,15,6)$ design has trivial automorphism group.

A surprising result of the switching of $W$ with respect to a row with weight structure $(15,21)$
that was suggested by the weight distribution in Table \ref{tab2} and was proved by inspection,
is the following.

\begin{thm}
\label{t1}

(a) The 36 codewords of $W_{408}$ with weight structure $(15,21)$ form a regular Hadamard matrix $H$ 
     which is monomially equivalent to the Paley-Hadamard matrix of type II. 
     
 (b) The symmetric 2-$(36,15,6)$ design $D$ associated with $H$ has a full automorphism
      group of order 24.
      
      (c) The incidence matrix of $D$ has 3-rank 18, and its linear span
      over $GF(3)$ is a code equivalent to the Pless symmetry code $C(17)$.
  \end{thm}
  
  \begin{note}
  
  {\rm
  Every row of the regular Hadamard matrix $H$ from Theorem \ref{t1} contains 15 entries equal to 1 and
  21 entries equal to $-1$. A $(0,1)$-incidence matrix $A$ of the associated symmetric 2-$(36,15,6)$ design
$D$ is obtained by  adding the all-one codeword $\bar{1}$ to every row of $H$,
followed by a multiplication of all rows by $2 \pmod 3$.  Hence, the ternary code
spanned by the rows of $H$ contains also the rows of $A$.
}
\end{note}

 A  $(0,1)$-incidence matrix of $D$  is given in the Appendix.
   
\section{Symmetric 2-$(36,15,6)$ designs with an involution and their ternary codes}

It was proved in \cite[Theorem 5.1]{EN} that if $C$ is an extremal ternary self-dual code
of length 36 then either $C$ is equivalent to the Pless symmetry code or the full automorphism group
of $C$ is a subgroup of the cyclic group of order 8. 
This result and the fact that
the Pless symmetry code is equivalent to codes spanned by the incidence matrices of two nonisomorphic 2-$(36,15,6)$ designs, motivated our study of symmetric  2-$(36,15,6)$ designs with 
an automorphism of order 2 (or, an {\it involution}), and the ternary codes spanned by their incidence matrices.

In what follows we summarize the classification of symmetric  2-$(36,15,6)$ designs 
that are invariant under an involution and the ternary  code of length 36 spanned by
their incidence matrix is a self-dual code with minimum weight 12.

It is  known that every finite group of order $2^n$ contains a
subgroup of order $2^i$ for every $i$ in the range $1 \le i \le n$
(cf. \cite[Theorem 6.5, page 116]{rose}).
By this property, finding all nonisomorphic $2$-$(36,15,6)$ designs which are invariant under an involution and whose incidence matrix spans a ternary extremal self-dual code, will also give the enumeration of all nonisomorphic $2$-$(36,15,6)$ designs invariant under a nontrivial subgroup of the cyclic group of order $8$ with this property.

For the construction of $2$-$(36,15,6)$ designs we use the method for constructing orbit matrices with presumed action of an automorphism group, which are then indexed to construct designs (see, e.g., \cite{c-r}, \cite{cep}). After constructing $2$-$(36,15,6)$ designs, we check the $3$-rank of their incidence matrices, and if it is equal to $18$ we determine the minimum weight of the corresponding ternary code. In our work we use computers. In addition to our own computer programs, we use computer programs by V. \'{C}epuli\'{c} 
for the construction of orbit matrices and the computer algebra system MAGMA \cite{magma} when working with codes.\\
 
The first step in the construction of $2$-$(36,15,6)$ designs that admit an involution is to determine the possible orbit lengths distributions. For that we use the following statements.

\begin{thm}\cite[Corollary~3.7]{l}\label{cor-fp}
Suppose that a nonidentity automorphism $\sigma$ of a symmetric $2$-$(v,k, \lambda)$ design fixes $f$ points. Then
$$f \le v-2(k-\lambda) \qquad {\rm and} \qquad f \le ( \frac{\lambda}{k- \sqrt{k-\lambda}} ) v.$$
\end{thm}

\begin{thm}\cite[Proposition~4.23]{l}\label{involution}
Suppose that ${\mathcal D}$ is a nontrivial symmetric $2$-$(v,k, \lambda)$ design, with an involution $\sigma$fixing $f$ points and blocks. If $f \neq 0$, then
$$
f \ge \left \{
\begin{tabular}{l  l}
 $1 + \frac{k}{\lambda}$,    & if $k$ and $\lambda$ are both even, \\
 $1 + \frac{k-1}{\lambda}$,  & otherwise. \\
\end{tabular} \right .
$$
\end{thm}

It follows that an involution acting on $2$-$(36,15,6)$ design could have $f$ fixed points, where $f \in \{4, 6, 8, 10, 12, 14, 16, 18 \}$.

Our analysis shows that orbit matrices do not exist for $f \in \{6, 14, 18 \}$. 
The results for the remaining cases  are summarized in Table \ref{tab3} and the constructed orbit matrices are available at
\begin{verbatim}
 https://www.math.uniri.hr/~sanjar/structures/
 \end{verbatim}
 The third row of Table \ref{tab3} contains information on the maximal minimum weight among the self-dual ternary codes spanned by incidence matrices of the corresponding $2$-$(36,15,6)$ designs, where symbol $\times$ indicates that all constructed designs have incidence matrices with 3-rank smaller than 18. \\

\begin{table}[htpb!]
\begin{center} \begin{footnotesize}
\begin{tabular}{|c ||c| c | c|c|c| }
 \hline 
 
Number of fixed points& 4 & 8 & 10 & 12 & 16 \\
 \hline
Number of orbit matrices& 12991& 670 & 56 & 311 & 83   \\
\hline
max $d$ in self-dual codes   & 12 & $\times$ & $\times$ & 9& $\times$ \\
\hline 
\end{tabular} \end{footnotesize}
 \caption{Symmetric 2-$(36,15,6)$ designs with an involution}\label{tab3}
\end{center} 
\end{table}
 
All $2$-$(36,15,6)$ designs whose incidence matrix span an extremal ternary code are mutually isomorphic. They are also isomorphic to the design discussed in Section 
\ref{sec2}, whose incidence matrix is given in the Appendix.

The results of our classification of symmetric 2-$(36,15,6)$ designs with the desired properties
can be summarized as follows.

\begin{thm}

(a) Up to isomorphism, there exists exactly one symmetric 2-$(36,15,6)$ design $D$
that admits an automorphism of order 2 and its incidence matrix spans
an extremal ternary self-dual code of length 36.

(b) The full automorphism group $G$ of $D$ is of order 24, 
and  $G$ is isomorphic to the symmetric group $S_4$.

(c)  The regular Hadamard matrix associated with $D$ is equivalent to 
       the Paley-Hadamard matrix of type II. 
       
(d) The ternary code spanned by $D$
is equivalent to the Pless symmetry code.

\end{thm}

\bigskip
\noindent {\bf Statements and Declarations} 

\noindent {\bf Competing Interests:} The authors declare no conflict of interest.\\
{\bf Contributions:} This is a joint collaboration with both authors contributing substantially throughout.\\
{\bf Data Availability:} 
The links to data generated and analysed are included in this article.\\
{\bf Funding:} The first author is supported by {\rm C}roatian Science Foundation under the project 6732.

\section{Appendix}

\begin{scriptsize}
\begin{verbatim}

 1 1 1 0 1 1 1 1 1 1 1 1 1 1 1 1 0 0 0 0 0 0 0 0 0 0 0 0 0 0 0 0 0 0 0 0
 1 1 0 1 1 1 1 1 0 0 0 0 0 0 0 0 1 1 1 1 1 1 1 1 0 0 0 0 0 0 0 0 0 0 0 0
 1 0 1 1 0 0 0 0 1 1 1 1 0 0 0 0 1 1 1 1 0 0 0 0 1 1 1 1 0 0 0 0 0 0 0 0
 0 1 1 1 0 0 0 0 0 0 0 0 1 1 1 1 0 0 0 0 1 1 1 1 1 1 1 1 0 0 0 0 0 0 0 0
 1 1 0 0 1 0 1 0 1 0 0 0 1 0 0 0 1 0 0 0 1 0 0 0 1 1 1 0 1 1 1 0 1 0 0 0
 1 1 0 0 0 1 0 1 0 1 0 0 0 1 0 0 0 1 0 0 0 1 0 0 1 1 0 1 1 1 0 1 0 1 0 0
 1 1 0 0 1 0 0 1 0 0 1 0 0 0 1 0 0 0 1 0 0 0 1 0 1 0 1 1 0 0 0 1 1 0 1 1
 1 1 0 0 0 1 1 0 0 0 0 1 0 0 0 1 0 0 0 1 0 0 0 1 0 1 1 1 0 0 1 0 0 1 1 1
 1 0 1 0 1 0 0 0 1 0 0 0 1 1 0 0 0 1 0 1 0 1 0 1 0 0 0 1 0 0 1 1 1 0 1 0
 1 0 1 0 0 1 0 0 0 1 0 0 1 1 0 0 1 0 1 0 1 0 1 0 0 0 1 0 0 0 1 1 0 1 0 1
 1 0 1 0 0 0 1 0 0 0 0 1 0 0 1 1 0 1 1 0 0 1 1 0 1 0 0 0 1 0 1 0 1 1 0 0
 1 0 1 0 0 0 0 1 0 0 1 0 0 0 1 1 1 0 0 1 1 0 0 1 0 1 0 0 0 1 0 1 1 1 0 0
 1 0 0 1 1 0 0 0 1 0 1 0 0 1 0 1 0 0 1 0 1 1 0 0 0 1 0 0 1 0 0 0 0 1 1 1
 1 0 0 1 0 1 0 0 0 1 0 1 1 0 1 0 0 0 0 1 1 1 0 0 1 0 0 0 0 1 0 0 1 0 1 1
 1 0 0 1 0 0 1 0 1 0 0 1 0 1 1 0 1 0 0 0 0 0 1 1 0 0 0 1 1 1 0 1 0 0 0 1
 1 0 0 1 0 0 0 1 0 1 1 0 1 0 0 1 0 1 0 0 0 0 1 1 0 0 1 0 1 1 1 0 0 0 1 0
 0 1 1 0 1 0 0 0 0 0 1 1 1 0 0 0 0 1 0 1 1 0 1 0 0 0 0 1 1 1 0 0 0 1 0 1
 0 1 1 0 0 1 0 0 0 0 1 1 0 1 0 0 1 0 1 0 0 1 0 1 0 0 1 0 1 1 0 0 1 0 1 0
 0 1 1 0 0 0 1 0 1 1 0 0 0 0 0 1 0 1 1 0 1 0 0 1 1 0 0 0 0 1 0 1 0 0 1 1
 0 1 1 0 0 0 0 1 1 1 0 0 0 0 1 0 1 0 0 1 0 1 1 0 0 1 0 0 1 0 1 0 0 0 1 1
 0 1 0 1 1 0 0 0 0 1 0 1 0 1 0 1 0 0 1 1 0 0 1 0 0 1 0 0 0 1 1 1 1 0 0 0
 0 1 0 1 0 1 0 0 1 0 1 0 1 0 1 0 0 0 1 1 0 0 0 1 1 0 0 0 1 0 1 1 0 1 0 0
 0 1 0 1 0 0 1 0 0 1 1 0 0 1 1 0 1 1 0 0 1 0 0 0 0 0 0 1 0 0 1 0 1 1 1 0
 0 1 0 1 0 0 0 1 1 0 0 1 1 0 0 1 1 1 0 0 0 1 0 0 0 0 1 0 0 0 0 1 1 1 0 1
 0 0 1 1 1 1 1 0 0 0 1 0 0 0 1 0 0 1 0 0 0 1 0 0 0 1 1 0 0 1 1 1 0 0 0 1
 0 0 1 1 1 1 0 1 0 0 0 1 0 0 0 1 1 0 0 0 1 0 0 0 1 0 0 1 1 0 1 1 0 0 1 0
 0 0 1 1 1 0 1 1 0 1 0 0 0 1 0 0 0 0 0 1 0 0 0 1 1 0 1 0 1 0 0 0 1 1 0 1
 0 0 1 1 0 1 1 1 1 0 0 0 1 0 0 0 0 0 1 0 0 0 1 0 0 1 0 1 0 1 0 0 1 1 1 0
 0 0 0 0 1 1 0 0 1 0 0 0 0 1 1 1 1 1 0 1 0 0 1 0 1 0 1 0 0 1 0 0 0 1 1 0
 0 0 0 0 1 1 0 0 0 1 0 0 1 0 1 1 1 1 1 0 0 0 0 1 0 1 0 1 1 0 0 0 1 0 0 1
 0 0 0 0 1 0 0 1 1 1 0 1 0 0 1 0 0 0 1 0 1 1 0 1 0 0 1 1 0 1 1 0 0 1 0 0
 0 0 0 0 0 1 1 0 1 1 1 0 0 0 0 1 0 0 0 1 1 1 1 0 0 0 1 1 1 0 0 1 1 0 0 0
 0 0 0 0 1 0 1 0 0 1 1 1 1 0 0 0 1 0 0 0 0 1 1 1 1 1 0 0 0 0 0 1 0 1 1 0
 0 0 0 0 0 1 0 1 1 0 1 1 0 1 0 0 0 1 0 0 1 0 1 1 1 1 0 0 0 0 1 0 1 0 0 1
 0 0 0 0 0 0 1 1 0 0 1 0 1 1 0 1 1 0 1 1 0 1 0 0 1 0 0 1 0 1 1 0 0 0 0 1
 0 0 0 0 0 0 1 1 0 0 0 1 1 1 1 0 0 1 1 1 1 0 0 0 0 1 1 0 1 0 0 1 0 0 1 0
\end{verbatim}
\end{scriptsize}

\begin{center}
A symmetric 2-$(36,15,6)$ design associated with the Paley-Hadamard matrix of type II
\end{center}

\end{document}